\documentclass[graybox]{svmult}

\usepackage{mathrsfs}    
\usepackage{stmaryrd}
\usepackage{eucal}     
\usepackage[all]{xy}

\usepackage{type1cm}        
\usepackage{makeidx}         
\usepackage{graphicx}        
\usepackage{multicol}        
\usepackage[bottom]{footmisc}

\usepackage{newtxtext}       %
\usepackage[varvw]{newtxmath}       

\usepackage{tikz}
\newcommand*\circled[1]{\tikz[baseline=(char.base)]{
            \node[shape=circle,draw,inner sep=0.05pt] (char) {#1};}}

\newcommand{\Z}{\mathbf{Z}}

\newcommand{\kk}{\mathbf{k}}
\newcommand{\sO}{\mathscr{O}}
\newcommand{\sI}{\mathscr{I}}
\def\nP{\mathbf{P}}
\def\id{\operatorname{id}}
\def\coker{\operatorname{coker}}

\def\reg{\operatorname{reg}}
\def\depth{\operatorname{depth}}
\def\pd{\operatorname{pd}}

\def\codim{\operatorname{codim}}

\def\Sing{\operatorname{Sing}}
\def\Supp{\operatorname{Supp}}

\DeclareMathOperator{\pr}{pr} 
\DeclareMathOperator{\Cliff}{Cliff} 
\DeclareMathOperator{\gon}{gon}

\makeindex             

\begin{document}

\title*{Syzygies of algebraic varieties through symmetric products of algebraic curves}
\author{Jinhyung Park\orcidID{0000-0002-9710-9728}}
\institute{Jinhyung Park \at Department of Mathematical Sciences, KAIST, 291 Daehak-ro, Yuseong-gu, Daejeon 34141, Replubic of Korea, \email{parkjh13@kaist.ac.kr}}
%
%
\maketitle

\abstract{This is a survey paper on recent work on syzygies of algebraic varieties. We discuss the gonality conjecture on weight-one syzygies of algebraic curves, syzygies of secant varieties of algebraic curves, syzygies of tangent developable surfaces and Green's conjecture on syzygies of canonical curves, and asymptotic syzygies of algebraic varieties. All results considered in this paper were proven using the geometry of symmetric products of algebraic curves.}

\section{Introduction}

Throughout the paper, we work over an algebraically closed field $\kk$ of arbitrary characteristic. 
In the late 1960s, Mumford \cite{Mumford} proved that if $X$ is a smooth projective variety and $L$ is a sufficiently positive line bundle, then $X \subseteq \nP H^0(X, L) = \nP^r$ is projectively normal and $I_{X|\nP^r}$ is generated by quadrics. In the early 1980s, Green realized that the classical results on equations defining algebraic varieties should be seen as the first cases of a much more general picture involving higher syzygies. After his pioneering work \cite{Green1, Green2} on Koszul cohomology, there has been a considerable amount of interest in syzygies of algebraic varieties. In the early 2000s, Voisin \cite{Voisin1, Voisin2} resolved Green's conjecture for general canonical curves using the Hilbert scheme interpretation of Koszul cohomology. Voisin's method has led to a great deal of progress on syzygies of algebraic varieties over the past two decades. Some of the recent developments were surveyed in \cite{EL3} and \cite{Farkas1}. We will discuss further advances centered around the author's work since those papers.

Let $X$ be a projective scheme of dimension $n$, and $L$ be a globally generated line bundle on $X$ with $r:=h^0(X, L)-1$. We usually assume that $L$ is sufficiently positive; in particular, $L$ is very ample, so it gives rise to an embedding $X \subseteq \nP H^0(X, L) = \nP^{r}$. Let $B$ be a coherent sheaf on $X$ with $H^0(X, B \otimes L^{-m})=0$ for $m > 0$. By Hilbert syzygy theorem, the finitely generated graded section module
$$
R=R(X, B; L):=\bigoplus_{m \in \Z} H^0(X, B \otimes L^m)
$$
over the polynomial ring $S:=\bigoplus_{m \in \Z} S^m H^0(X, L)$ has a minimal free resolution
$$
0 \longleftarrow R \longleftarrow E_0 \longleftarrow E_1 \longleftarrow \cdots \longleftarrow E_r \longleftarrow 0,
$$
where
$$
E_p = \bigoplus_{q \in \Z} \underbrace{K_{p,q}(X, B;L)}_{\emph{Koszul cohomology}}  \otimes_{\kk} S(-p-q).
$$
Here $K_{p,q}(X,B;L)$ is the cohomology of the Koszul-type complex
$$
\wedge^{p+1}H^0(L)\otimes H^0(B\otimes L^{q+1}) \! \rightarrow \!  \wedge^pH^0(L)\otimes H^0(B\otimes L^q) \!  \rightarrow \!  \wedge^{p-1}H^0(L)\otimes H^0(B\otimes L^{q-1}),
$$
and it can be regarded as the space of $p$-th syzygies of weight $q$. Denote by
$$
\kappa_{p,q}=\kappa_{p,q}(X, B; L):=\dim K_{p,q}(X, B; L)
$$
the \emph{graded Betti numbers} of $R=R(X,B;L)$. The \emph{Betti table} of $R$ is the following: 
\begin{center}
\begin{tabular}{c|cccc}
& $0$ & $1$ & $2$ & $~\cdots~$ \\ \hline
~$0$~ & $\kappa_{0,0}$ & $\kappa_{1,0}$ & $\kappa_{2,0}$ & $~\cdots~$  \\
~$1$~ & $\kappa_{0,1}$ & $\kappa_{1,1}$ & $\kappa_{2,1}$& $~\cdots~$  \\
~$2$~ & $\kappa_{0,2}$ & $\kappa_{1,2}$ & $\kappa_{2,2}$& $~\cdots~$  \\
$\vdots$ & $\vdots$ &  $\vdots$ & $\vdots$ &   \\ 
\end{tabular}
\end{center}
The width of the Betti table of $R$ is the \emph{projective dimension} $\pd(R)$. By Auslander--Buchsbaum formula, $\pd(R)=r-\depth(R)$. The height of the Betti table of $R$ is the \emph{Castelnuovo--Mumford regularity} $\reg(R)$ (see \cite[Subsection 1.8]{positivity}).

Our main concern is the overall shape of the Betti table of $R(X,B;L)$ when $L$ is sufficiently positive. Specifically, we investigate from a geometric perspective when $\kappa_{p,q}(X,B;L)$ are vanishing and what their values are if they are not vanishing. In this paper, we always assume that $H^i(X, B \otimes L^m)=0$ for $i>0$ and $m>0$. Then $\reg(R) \leq n+1$, so $K_{p,q}(X, B; L) = 0$ for $q \geq n+2$. It is well-known that
$$
K_{p,q}(X, B; L)=H^{q-1}(X, \wedge^{p+q-1}M_L  \otimes B \otimes L)~~\text{ for $p \geq 0$ and $q \geq 2$},
$$
where $M_L$ is the kernel bundle of the evaluation map $H^0(X, L)\otimes \sO_X \to L$. If $H^i(X, B \otimes L^m)=0$ for $1 \leq i \leq n-1$ and $m \in \Z$ (i.e., $R(X, B; L)$ is Cohen--Macaulay), then $K_{p,q}(X, B; L)=H^{q}(X, \wedge^{p+q}M_L  \otimes B)$ and $K_{p,q}(X, B; L) = K_{r-p-n, n+1-q}(X, \omega_X \otimes B^{-1}; L)^{\vee}$ for $p \geq 0$ and $q \geq 0$. The latter can be easily checked using Serre duality.
For more details, we refer to \cite[Chapter 2]{AN}, \cite[Section 3]{EL1}, \cite[Section 2]{Park3}.

For a moment, suppose that $B=\sO_X$. We put $R(X;L):=R(X,\sO_X;L)$ and $K_{p,q}(X; L):=K_{p,q}(X, \sO_X; L)$. We know $K_{0,0}(X; L)= \kk$ and $K_{p,0}(X;L)=0$ for $p \geq 1$. For integers $d \geq 2$ and $m \geq 0$, we say that $L$ satisfies \emph{$N_{d,m}$-property} (when $d=2$, we simply say  \emph{$N_m$-property}) if $K_{p,q}(X; L) = 0$ for $0 \leq p \leq m$ and $q \geq d$. This notion was introduced in \cite{EGHP} based on \cite{GL1, GL2} (our definition is slightly different but the same when $X \subseteq \nP^r$ is arithmetically normal, i.e., $R=S/I_{X|\nP^r}$). As $K_{0,1}(X;L)=0$, we see that $R=S/I_{X|\nP^r}$ if and only if $L$ satisfies $N_0$-property. In this case, $\kappa_{1,d}(X;L)$ is the number of minimal generators of $I_{X|\nP^r}$ of degree $d+1$, so $L$ satisfies $N_1$-property if and only if $I_{X|\nP^r}$ is generated by quadrics. 

Generalizing classical results of Castelnuovo, Mattuck, Mumford, Saint-Donat, Fujita on equations defining algebraic curves, Green \cite[Theorem 4.a.1]{Green1} proved that if $C$ is a smooth projective curve of genus $g$ and $\deg L \geq 2g+1+p$, then $L$ satisfies $N_p$-property. Green--Lazarsfeld \cite[Theorem 1.2]{GL2} then showed that $K_{p+1,2}(C; L) \neq 0$ when $\deg L = 2g+1+p$ and $H^0(C, L \otimes \omega_C^{-1}) \neq 0$. Thus we obtain the following:

\begin{theorem}[{$\operatorname{char}(\kk)=0$}]\label{thm:Green}
If $\deg L \geq 2g+1$ and $H^0(C, L \otimes \omega_C^{-1}) \neq 0$, then 
$$
K_{p,2}(C; L) \neq 0~~\Longleftrightarrow~~r-g \leq p \leq r-1.
$$
\end{theorem}

\noindent The vanishing and nonvanishing of $K_{p,1}(C;L)$ is predicted by Green--Lazarsfeld's gonality conjecture \cite[Conjecture 3.7]{GL1}. This was verified by Ein--Lazarsfeld \cite{EL2} when $\deg L \gg 0$. Recently, Niu--Park \cite{NP} showed the following effective gonality theorem establishing a sharp effective bound on $L$ to satisfies the gonality conjecture.

\begin{theorem}[{$\operatorname{char}(\kk)=0$}]\label{thm:effgon}
If $\deg L \geq 2g+\gon(C)$, then 
$$
K_{p,1}(C; L) \neq 0 ~~\Longleftrightarrow~~ 1 \leq p \leq r-\gon(C).
$$
\end{theorem}

\noindent Here $\gon(C)$ is the \emph{gonality} of $C$, the minimal degree of a branch covering $C \to \nP^1$. As Castryck \cite{Castryck} observed that the gonality conjecture may not hold for a plane curve when $\deg L = 2g+\gon(C)-1$, this theorem is the best possible and thus gives a complete answer to the effective gonality problem. By Theorems \ref{thm:Green} and \ref{thm:effgon}, if $\deg L \geq 3g$, then the Betti table of $R(C;L)$ is the following:
$$
\begin{array}{c|ccccccccccc}
& ~0~ & ~1~ & ~2~ & ~\cdots~ & r-g-1 & r-g & \cdots & r-\gon(X) & r+1-\gon(X) & \cdots &   r-1 \\ \hline
~0~ &1 & - & - & \cdots & - & - & \cdots & - & - & \cdots & -\\
~1~ & - & * & * & \cdots & * & * & \cdots & * & - & \cdots & - \\
~2~ & - & - & - & \cdots & - & * & \cdots & * & * & \cdots & *
\end{array}
$$
Here ``$*$'' indicates nonzero and ``$-$'' indicates zero. The number of nonzero entries in the last row is $g$, and the number of last zero entries in the row above is $\gon(C)-1$. This reveals the surprising connection between the intrinsic geometric properties of the algebraic curve $C$ and the extrinsic algebraic properties of the embedded curve $C \subseteq \nP^r$. The syzygies of algebraic curves will be discussed in Section \ref{sec:gonality}. 

Sidman--Vermeire \cite{SV} proposed a generalization of Green's $(2g+1+p)$-theorem to secant varieties. Suppose that $\deg L \geq 2g+2k+1$. The \emph{$k$-th secant variety} $\Sigma_k$ of $C \subseteq \nP^r$ is the union of $(k+1)$-secant $k$-planes to $C$ in $\nP^r$. Note that $\Sigma_0=C$. We have $\dim \Sigma_k = 2k+1$ and $\Sing \Sigma_k = \Sigma_{k-1}$ except when $C=\nP^1$ and $\deg L = 2k+1$ (in this case $\Sigma_k= \nP^{2k+1}$).  Ein--Niu--Park \cite{ENP1} proved the following theorem resolving the Sidman--Vermeire conjecture \cite[Conjecture 1.3]{SV} (the case of $k=0$ is Green's $(2g+1+p)$-theorem, and the case of $k=1$ was proved by Chou--Song \cite{CS}, Sidman--Vermeire \cite{SV}, Ullery \cite{Ullery}, Vermeire \cite{Vermeire}, Yang \cite{RuijieYang}).

\begin{theorem}[{$\operatorname{char}(\kk)=0$}]\label{thm:ENP}
The $k$-th secant variety $\Sigma_k$ has normal Du Bois singularities, $\Sigma_k \subseteq \nP^r$ is arithmetically Cohen--Macaulay, and $\sO_{\Sigma_k}(1)$ satisfies $N_{k+2, p}$-property as soon as $\deg L \geq 2g+2k+1+p$.
\end{theorem}

\noindent If $e:=\codim \Sigma_k = r-2k-1$, then the Betti table of $R(\Sigma_k; \sO_{\Sigma_k}(1))$ is the following:
$$
\begin{array}{c|ccccccccc}
& 0 & 1 & 2 & ~\cdots &e-g-1 & e-g &  e-g+1 & \cdots & ~e \\ \hline
0 & 1 & - & - & ~\cdots  & - & - & - & \cdots & ~-  \\
1 & - & -  & -  & ~\cdots & - & - & - & \cdots & ~- \\
\vdots & \vdots  & \vdots  & \vdots & & \vdots   & \vdots  & \vdots & & ~\vdots \\
k & - & - & - & ~\cdots & - & - & - & \cdots & ~- \\
k+1 & - & * & * & ~\cdots & * & * & ? & \cdots & ? \\
k+2 & - & - & - & ~\cdots & - & -  & ? & \cdots & ? \\
\vdots & \vdots  & \vdots  & \vdots & & \vdots  & \vdots &\vdots  & & ~\vdots\\
2k+2 & - & - & - & ~\cdots & -  & - & ? & \cdots & ?
\end{array}
$$
When $\deg L \gg 0$, Choe--Kwak--Park \cite{CKP} completely determined the remaining part of the Betti table, which is a natural generalization of the gonality conjecture.

\begin{theorem}[{$\operatorname{char}(\kk)=0$}]\label{thm:CKP}
For each $k+1 \leq q \leq 2k+2$, if $e-g+1 \leq p \leq e$, then 
$$
K_{p,q}(\Sigma_k; \sO_{\Sigma_k}(1)) \neq 0~~\Longleftrightarrow~~e - g+1 \leq p \leq e - \gamma^{2k+2-q}(C).
$$
\end{theorem}

\noindent Here $(\gamma^0(C)+0, \gamma^1(C)+1, \gamma^2(C)+2, \ldots)$ is the \emph{gonality sequence} of $C$ with $\gamma^1(C)+1=\gon(C)$. This was originally conjectured in \cite{CK}. Theorem \ref{thm:CKP} says that the part marked with ``$?$'' in the above Betti table is the following:
\begin{center}
\begin{tabular}{c}
 \texttt{$*$ $\cdots$ $*$} $\overbrace{\text{$\, -$ $\cdots$ $-$ } \, \text{ $-$ $\, \cdots \, $ $-$ } \, \text{ $-$ $\, \cdots \,$ $-\,$}}^{\gamma^{k+1}(C)}$\\
\vdots\\
 \texttt{$*$ $\cdots$ $*$ $*$ $\cdots$ $*$}  $\overbrace{\text{$\, -$ $\, \cdots\, $ $-$ } \, \text{ $-$ $\, \cdots \,$ $-\,$}}^{\gamma^2(C)}$\\
 \texttt{$*$ $\cdots$ $*$ $*$ $\cdots$ $*$ $*$ $\cdots$ $*$} $\overbrace{\text{$\,-$ $ \, \cdots \,$ $-\,$}}^{\gamma^1(C)}$\\
$\underbrace{\texttt{$*$ $\cdots$ $*$ $*$ $\cdots$ $*$ $*$ $\cdots$ $*$ $*$ $\cdots$ $*$}}_{g}$\\
\end{tabular}
\end{center}
The syzygies of secant varieties will be discussed in Section \ref{sec:secant}.

The \emph{tangent developable surface} of $C \subseteq \nP^r$ is the union of tangent lines to $C$ in $\nP^r$. There has been a lot of interest in the syzygies of the tangent developable surface of a rational normal curve of degree $g$ because they are the same to the syzygies of general canonical curves of genus $g$. This was independently observed by O’Grady and Buchweitz–Schreyer in 1980s. After 35 years, Aprodu--Farkas--Papadima--Raicu--Weyman \cite{AFPRW} finally proved the following:

\begin{theorem}[{$\operatorname{char}(\kk) = 0$ \normalfont{or} $\geq (g+2)/2$}]\label{thm:AFPRW}
If $T$ is the tangent developable surface of a rational normal curve of degree $g$, then $\sO_T(1)$ satisfies $N_{\lfloor (g-3)/2 \rfloor}$-property.
\end{theorem}

\noindent This implies generic Green's conjecture \cite[Conjecture 5.6]{Green1}, which asserts that the canonical line bundle of a general curve of genus $g$ satisfies $N_{\lfloor (g-3)/2 \rfloor}$-property. Recently, it was realized by Park  \cite{Park2} that $T$ is a Weil divisor linearly equivalent to $-K_{\Sigma}$, where $\Sigma$ is the first secant variety of the rational normal curve of degree $g$. This leads to a substantialy simple geometric proof of Theorem \ref{thm:AFPRW}. Furthermore, it was also shown that if $C$ is a smooth projective curve of genus $g \geq 1$ and $L$ is a line bundle on $C$ with $\deg L \geq 4g+3$, then the tangent developable surface $T$ of $C$ in $\nP H^0(C, L)$ is arithmetically normal (\cite[Theorem 1.3]{Park2}). The syzygies of tangent developable surfaces and Green's conjecture will be discussed in Section \ref{sec:tangent}.

Green's $(2g+1+p)$-theorem has stimulated further work on $N_p$-property, and several analogous statements for higher dimensional algebraic varieties have been established. On the other hand, in 2001, Ottaviani--Paoletti \cite{OP} called attention to the failure of $N_m$-property. Putting $r_d:=h^0(\nP^2, \sO_{\nP^2}(d))-1$, they proved that
$$
K_{p,2}(\nP^2; \sO_{\nP^2}(d)) \neq 0~~\text{ for $3d-2 \leq p \leq r_d-3$}.
$$
As $r_d \approx d^2/2$, this shows $N_m$-property for $\sO_{\nP^2}(d)$ describes only a small fraction of the syzygies of the $d$-th Veronese embedding of $\nP^2$. It is an interesting problem to describe the overall asymptotic behaviors of $K_{p,q}(X, B; L)$ as the positivity of $L$ grows. To set the stage, assume that $X$ is a smooth projective variety of dimension $n$ and $B$ is a line bundle on $X$, and let $L_d:=\sO_X(dA+P)$ for an integer $d \gg 0$ and $r_d:=h^0(X, L_d)-1$, where $A$ is an ample divisor on $X$ and $P$ is an arbitrary divisor on $X$. By \cite[Proposition 5.1 and Corollary 5.2]{EL2}, we have
\begin{equation}\label{eq:q=0,n+1}
\begin{array}{l}
K_{p,0}(X, B; L_d) \neq 0 ~\Longleftrightarrow~0 \leq p \leq h^0(X, B)-1;\\
K_{p, n+1}(X, B; L_d) \neq 0 ~\Longleftrightarrow~r_d-n-h^0(X, \omega_X \otimes B^{-1}) + 1 \leq p \leq r_d-n.
\end{array}
\end{equation}
The main issue is then to study vanishing and nonvanishing of $K_{p,q}(X, B; L_d)$ for $1 \leq q \leq n$. Ein--Lazarsfeld \cite{EL2} and Park \cite{Park1, Park3} established the following:

\begin{theorem}[{$\operatorname{char}(\kk)\geq 0$}]\label{thm:asysyz}
For each $1 \leq q \leq n$, there exist functions $c_q(d)$ and $c_q'(d)$ depending on $q, d$ and $X, B, A, P$ with
$$
c_q(d) = \Theta(d^{q-1})~~\text{ and }~~c_q'(d) = \begin{cases} \Theta(d^{n-q}) & \text{if $H^{q-1}(X, B) =0$ or $q=1$} \\  q-1 & \text{if $H^{q-1}(X, B) \neq 0$ and $q \geq 2$} \end{cases}
$$
such that if $d$ is sufficiently large, then
$$
K_{p,q}(X, B; L_d) \neq 0~~\Longleftrightarrow~~c_q(d) \leq p \leq r_d-c_q'(d).
$$
\end{theorem}

\noindent For a nonnegative function $f(d)$, we write $f(d) = \Theta(d^k)$ if there are constants $C_1, C_2 > 0$ such that $C_1 d^k \leq f(d) \leq C_2 d^k$ for any sufficiently large positive integer $d$ (for convention, we write $f(d) = \Theta(1)$ even when $f(d)=0$ for all $d \gg 0$). A  widely-accepted intuition derived from the case of curves was that syzygies become simpler as the positivity of the embedding line bundle is increasing. However, as was noticed by Ein--Lazarsfeld \cite{EL2}, this intuition was misleading since Theorem \ref{thm:asysyz} says that almost all possibly nonzero syzygies are actually nonzero. Instead, one may say that the syzygies have a surprisingly uniform asymptotic behavior as the positivity of the embedding line bundle is increasing. Although Theorem \ref{thm:asysyz} is for higher dimensional algebraic varieties, the geometry of symmetric products of algebraic curves plays a crucial role in an important step of the proof. The asymptotic syzygies of algebraic varieties will be discussed in Section \ref{sec:asymptotic}.

\section{Gonality conjecture on weight-one syzygies of algebraic curves}\label{sec:gonality}
In this section, we assume that $\operatorname{char}(\kk) = 0$.
Let $C$ be a smooth projective curve of genus $g \geq 0$, and $B,L$ be line bundles on $C$. We begin with the preparation for the necessary techniques relating $K_{p,1}(C,B;L)$ to vector bundles on symmetric products of $C$. For more details, we refer the reader to  \cite{EL2}, \cite{ENP1}, \cite{NP}. 

For an integer $m \geq 0$, the $m$-th symmetric product $C_m$ of $C$ parameterizes effective divisors of degree $m$ on $C$. It is well-known that $C_m$ is a smooth projective variety of dimension $m$. For integers $m,n \geq 1$, there is an addition map
$$
\sigma_{m,n} \colon C_{m}\times C_n \longrightarrow C_{m+n},~~(\xi,\xi') \longmapsto \xi+\xi',
$$
which is a finite morphism.
The incidence subvariety 
$$
D_{m,n} = \{(\xi, \eta) \in C_{m} \times C_{n} \mid \Supp(\xi) \cap \Supp(\eta) \neq \emptyset\} \subseteq C_m \times C_n
$$ 
is an effective divisor on $C_m \times C_n$. If $n=1$, then $D_{m,1}=C_{m-1} \times C$ (the map is given by $(\xi,x) \mapsto (\xi-x, x)$) and $\sigma_{m-1,1}  = \pr_1|_{D_{m,1}} \colon C_{m-1} \times C \to C_m$ is the universal family over $C_{m}$, where $\pr_1 \colon C_m \times C \to C_m$ is the projection. Let 
$$
E_{m,L}:=\sigma_{m-1,1,*} (\sO_{C_{m-1}} \boxtimes L)
$$
be the  \emph{tautological bundle}, which is a vector bundle of rank $m$ on $C_{m}$ with $H^0(C_m, E_{m,L})=H^0(C, L)$. Put $N_{m,L}:=\det (E_{m,L})$. If we write $N_{m,\sO_C}=\sO_{C_{m}}(-\delta_{m})$, then $\sigma_{m-1,1}^*\delta_{m} = \pr_1^* \delta_{m-1} +
D_{m-1,1} \subseteq C_{m-1} \times C$. Note that $\Omega_{C_m}^1=E_{m, \omega_C}$ and $\omega_{C_m} = N_{m, \omega_C}$. Now, consider the group action of the symmetric group $\mathfrak{S}_{m}$ on the $m$-th ordinary product $C^m$ of $C$ and the quotient map $q_m \colon C^m \to C_m$. There is a line bundle $S_{m, L}$ on $C_{m}$ which is the $\mathfrak{S}_{m}$-invariant descent of 
$$
L^{\boxtimes m}:=\underbrace{L\boxtimes  \cdots \boxtimes L}_{\text{$m$ times}}.
$$ 
Note that $q_m^* S_{m,L} = L^{\boxtimes m}$ and $N_{m, L} = S_{m, L}(-\delta_{m})$. By \cite[Lemma 3.7]{ENP1} (see also \cite[Lemma 2.4]{Agostini2}), we have
$$
\begin{array}{l}
H^i(C_m, N_{m,L}) = \wedge^{m-i} H^0(C, L) \otimes S^i H^1(C, L);\\
H^i(C_m, S_{m,L}) = S^{m-i} H^0(C, L) \otimes \wedge^i H^1(C, L).
\end{array}
$$
Notice that the fiber of $E_{p+1, B}$ over $\xi \in C_{p+1}$ is $H^0(\xi, B|_{\xi})$. Thus $B$ is $p$-very ample (i.e., the restriction map $H^0(C, B) \to H^0(\xi, B|_{\xi})$ is surjective for every $\xi \in C_{p+1}$) if and only if $\operatorname{ev}_{p+1, B}$ is surjective if and only if $E_{p+1, B}$ is globally generated. From now on, suppose that $B$ is $p$-very ample and $L$ is globally generated. Let $M_{p+1, B}$ be the kernel bundle of the evaluation map $H^0(C, B) \otimes \sO_{C_{p+1}} \to E_{p+1, B}$. The key observation of Voisin (see e.g., \cite[Lemma 1.1]{EL2}) is that $K_{p,1}(C, B; L)$ is the cokernel of the multiplication map
$$
H^0(C, B)\otimes H^0(C_{p+1}, N_{p+1,L})\longrightarrow H^0(C_{p+1}, E_{p+1,B}\otimes N_{p+1,L}).
$$
Thus $H^1(C_{p+1}, M_{p+1, B} \otimes N_{p+1, L})=0$ implies $K_{p,1}(C, B; L) = 0$, and the converse also holds when $H^1(C,L)=0$.

We now turn to the gonality conjecture of Green--Lazarsfeld \cite[Conjecture 3.7]{GL1} asserting that $K_{p,1}(C;L) = 0$ for $p \geq r-\gon(C)+1$ when $\deg L \gg 0$.
By the duality theorem (\cite[Theorem 2.c.1]{Green1}), 
$$
K_{p,1}(C; L)=K_{r(L)-p-1, 1}(C, \omega_C; L)^{\vee},
$$ 
where $r(L):=h^0(C, L)-1$. It is well-known that $\gon(C) - 2 \geq p$ if and only if $\omega_C$ is $p$-very ample. Thus the gonality conjecture can be restated as if $\omega_C$ is $p$-very ample, then $K_{p,1}(C, \omega_C; L)=0$. Niu--Park \cite[Theorem 1.3]{NP} established a general effective vanishing theorem for weight-one syzygies of algebraic curves.

\begin{theorem}\label{thm:NP}
Suppose that $B$ is $p$-very ample. If 
\begin{equation}\label{eq:effvansyzcur}
h^1(C, L \otimes B^{-1}) \leq r(B)-p-1,
\end{equation}
then $L$ is $p$-very ample and $K_{p,1}(C, B; L)=0$.
\end{theorem}

This theorem is a significant improvement of Rathmann's theorem \cite[Theorem 1.2]{Rathmann1} saying if $h^1(C, L) = h^1(C, L \otimes B^{-1})=0$, then $K_{p,1}(C, B; L)=0$. This implies that if $\deg L \geq 4g-3$, then the gonality conjecture holds. By Riemann--Roch, the condition (\ref{eq:effvansyzcur}) is equivalent to 
$$
\deg L \geq 2g+p+h^0(C, L \otimes B^{-1})-h^1(C, B).
$$
It is worth noting that Theorem \ref{thm:NP} immediately implies Green's $(2g+1+p)$-theorem since $K_{p,2}(C; L) = K_{p,1}(C, L; L)$. Using Theorem \ref{thm:NP} and Green--Lazarsfeld's nonvanishing theorem \cite[Appendix]{Green1} (see also Proposition \ref{prop:effnonvanK_{p,1}(C,B;L)}), Niu--Park proved the following effective gonality theorem \cite[Corollary 5.3]{NP}, which implies Theorem \ref{thm:effgon} saying that the gonality conjecture holds when $\deg L \geq 2g+\gon(C)$.

\begin{theorem}\label{thm:effgondualform}
Assume that $g \geq 2$ and $\deg L \geq 2g+p$. Then $K_{p,1}(C, \omega_C; L) \neq 0$ if and only if one of the following holds:
\begin{enumerate}
	\item[$(a)$] $\omega_C$ is not $p$-very ample.
	\item[$(b)$] $C \subseteq  \nP H^0(C, H) = \nP^2$ is a plane curve of degree $p+3 \geq 4$ and $L=\omega_C \otimes H$. In this case, $\omega_C$ is $p$-very ample and $\deg L = 2g+p+1$.
	\item[$(c)$] $C$ is arbitrary and $L = \omega_C(\xi)$ for some $\xi \in C_{p+2}$ with $\dim |\xi| = 1$. In this case, $\omega_C$ is $p$-very ample and $\deg L = 2g+p$.
\end{enumerate}
\end{theorem}

It was originally expected by Green--Lazarsfeld \cite[page 87]{GL1} and explicitly conjectured by Farkas--Kemeny \cite[page 3]{FK2} that $K_{p,1}(C, \omega_C; L)\neq 0$ if and only if $\omega_C$ is not $p$-very ample when $\deg L \geq 2g+p+1$. However, Castryck \cite{Castryck} noticed the example in $(b)$ and suggested that this conjecture should be modified. Theorem \ref{thm:effgondualform} completely resolves the effective gonality problem.

To prove Theorem \ref{thm:NP}, we need to show that $H^1(C_{p+1}, M_{p+1,B} \otimes N_{p+1, L})=0$. By Serre duality, it is equivalent to
$$
H^p(C_{p+1}, \wedge^{r(B)-p-1} M_{p+1, B} \otimes N_{p+1, \omega_C \otimes B \otimes L^{-1}})=0.
$$
By \cite[Lemma 2.3]{NP} (cf. \cite[Lemma 3.3]{CKP}), we have
\begin{align*}
R^i \pr_{2,*} (N_{r(B)-p-1, B} \boxtimes \sO_{C_{p+1}})(-D_{r(B)-p-1, p+1})=S^i H^1(B) \otimes  \wedge^{r(B)-p-1-i} M_{p+2, B}
\end{align*}
for $i \geq 0$, where $\pr_2 \colon  C_{r(B)-p-1} \times C_{p+1} \to C_{p+1}$ is the projection map. By the Leary spectral sequence for $\pr_2$,  the desired cohomology vanishing is implied by 
$$
H^i(C_j \times C_{p+1}, (N_{j, B} \boxtimes N_{p+1, \omega_C \otimes B \otimes L^{-1}})(-D_{r(B)-p-1, p+1})) = 0
$$
for $0 \leq j \leq r(B)-p-1$ and $0 \leq i \leq p$. Taking the dual and considering the Leray spectral sequence for the projection $\pr_1 \colon C_j \times C_{p+1} \to C_j$, we see that the desired cohomology vanishing is implied by
$$
\dim \mathcal{L}_j^{\ell}(\omega_C \otimes B \otimes L^{-1})  \leq p+j-i-\ell~~\text{ for $\ell \geq 1$}.
$$
But this can be easily deduced from the condition (\ref{eq:effvansyzcur}). Here 
$$
\mathcal{L}_j^{\ell}(\omega_C \otimes B \otimes L^{-1}) := \{\xi \in C_j \mid h^1(C, (L \otimes B^{-1})(-\xi)) \geq \ell\},
$$ 
which is the \emph{secant space} of $j$-secant $(r(\omega_C \otimes B \otimes L^{-1})-\ell)$-planes to $C$ in $\nP H^0(C, \omega_C \otimes B \otimes L^{-1})$ when $\ell \geq 1$ and $j \geq r(\omega_C \otimes B \otimes L^{-1})-\ell+1$. 

Niu--Park \cite[Corollary 4.7]{NP} also showed a general effective nonvanishing theorem for weight-one syzygies of algebraic curves (see also \cite[Proposition 3.5]{Agostini1}).

\begin{proposition}\label{prop:effnonvanK_{p,1}(C,B;L)}
Suppose that $B$ is $p$-very ample but not $(p+1)$-very ample. 
If $\deg L\geq 2g+p+1$ then $K_{p+1,1}(C, B;L) \neq 0$ except for the following:
\begin{enumerate}
\item[$(a)$] $g=0$, $B=\sO_{\nP^1}(p)$, $L=\sO_{\nP^1}(p+1)$ with $p \geq 1$. In this case, $K_{p+1,1}(\nP^1,B; L)=0$.
\item[$(b)$] $g=1$, $B=\sO_C$, $\deg L=3$. In this case, $K_{1,1}(C; L)=0$.
\end{enumerate}
\end{proposition}

Finally, we discuss Green--Lazarsfeld's secant conjecture \cite[Conjecture 3.4]{GL1}: if $L$ is very ample and
$$
\deg L \geq 2g+1+p - 2h^1(C, L) - \Cliff(C),
$$
then $L$ satisfies $N_p$-property if and only if $L$ is $(p+1)$-very ample. This was verified by Green--Lazarsfeld for $p=0$ \cite[Theorem 1]{GL1}. Here $\Cliff(C)$ is the \emph{Clifford index} of $C$ defined by the minimum of  $\Cliff(A):=\deg A - 2h^0(C, A)+2$ for a line bundle $A$ on $C$ with $h^0(C, A) \geq 2$ and $h^1(C, A) \geq 2$. In the situation of the conjecture, $h^1(C, L) \leq 1$. The conjecture for $h^1(C, L)=1$ is equivalent to Green's conjecture \cite[Conjecture 5.1]{Green1} (see \cite[Proposition 4.30]{AN}). When $h^1(C, L)=0$, the conjecture can be restated as if $\deg L =2g+1+p-c$ with $0 \leq c \leq \Cliff(C)$ and $L$ is $(p+1)$-very ample, then $L$ satisfies $N_p$-property (indeed, the converse direction is known (see \cite[Theorem 4.37]{AN})). This was shown by Farkas--Kemeny \cite{FK1} when $C$ and $L$ are general. The case of $c=0$ is Green's $(2g+1+p)$-theorem,  the case of $c=1$ was verified by Green--Lazarsfeld \cite[Theorem 1.2]{GL2}, and the case of $c=2$ was recently settled by Agostini \cite{Agostini2} when $C$ is not bielliptic. 

We briefly explain Agostini's approach in \cite{Agostini2}. We want to prove that 
$$
K_{p,2}(C;L)=H^1(C, \wedge^{p+1} M_L \otimes L) = H^0(C, \wedge^{g-c} M_L \otimes \omega_C)^{\vee} = 0.
$$
As $H^0(C, \wedge^{g-c} M_L \otimes \omega_C) = H^0(C_{g-c} \times C, (N_{g-c,L} \boxtimes \omega_C)(-D_{g-c,1}))$, it suffices to show that
$$
H^0(C_{g-2c} \times C_c \times C, (\sigma_{g-2c, c} \times \id_C)^* (N_{g-c,L} \boxtimes \omega_C)(-D_{g-c,1}))=0.
$$
This can be deduced from 
$$
H^0(C_c \times C, N_{c, L(-\xi)} \boxtimes \omega_C(-\xi))(-D_{c,1}))=0~~\text{ for general $\xi \in C_{g-2c}$}.
$$
When $c=0$, this is $H^0(C, \omega_C(-\xi))=0$ for general $\xi \in C_g$ but it is clear. Assume that $c \geq 1$ and $\dim W_{g-c}^1(C) = g-2c-2$. Agostini \cite[Lemma 4.2]{Agostini2} observed $\omega_C(-\xi)$ is $(c-1)$-very ample under this assumption. The above vanishing is equivalent to
$$
H^0(C_c, M_{c, \omega_C(-\xi)} \otimes N_{c, L(-\xi)})=0.
$$
As $M_{c, \omega_C(-\xi)} \otimes N_{c, \omega_C(-\xi)}=\wedge^{c-1} M_{c, \omega_C(-\xi)}^{\vee}$, we have a long exact sequence
$$
\begin{array}{l}
\cdots \longrightarrow \wedge^{c-3} H^0(C, \omega_C(-\xi)) \otimes S^2 E_{c, \omega_C(-\xi)}^{\vee}   \longrightarrow \wedge^{c-2} H^0(C, \omega_C(-\xi)) \otimes E_{c, \omega_C(-\xi)}^{\vee} \\
\longrightarrow  \wedge^{c-1} H^0(C, \omega_C(-\xi))^{\vee} \otimes \sO_C \longrightarrow M_{c, \omega_C(-\xi)} \otimes N_{c, \omega_C(-\xi)} \longrightarrow 0.
\end{array}
$$
Thus we reduce the problem to
\begin{equation}\label{eq:vansymS_E}
H^i(C_c, S^i E_{c, \omega_C(-\xi)}^{\vee} \otimes S_{c, L \otimes \omega_C^{-1}})=0~~\text{ for $0 \leq i \leq c-1$}.
\end{equation}
When $c=1$, this is $H^0(C, L \otimes \omega_C^{-1})=0$, which holds since $\deg L = 2g+p$ but $L$ is $(p+1)$-very ample. When $c=2$, Agostini proved this vanishing holds if and only if $g \neq p+4$ in \cite[Lemma 4.5]{Agostini2}. In the exceptional case of $g=p+4$, some additional arguments are required to settle Green--Lazarsfeld's secant conjecture (see \cite[Section 5]{Agostini2}). Recently, Farkas \cite{Farkas2} has given a different account of the secant conjecture for the case of $g=p+2+c$. It would be very nice if one can verify the vanishing (\ref{eq:vansymS_E}) for $c \geq 3$ when $g \neq p+2+c$.

\section{Syzygies of secant varieties of smooth projective curves}\label{sec:secant}
In this section, we assume that $\operatorname{char}(\kk) = 0$. Let $C$ be a smooth projective curve of genus $g \geq 0$, and $L$ be a line bundle on $C$ with $\deg L \geq 2g+2k+1$ for an integer $k \geq 0$. The tautological bundle $E_{k+1, L}$ on $C_{k+1}$ is globally generated with $H^0(C_{k+1}, E_{k+1,L}) = H^0(C,L)$, so $|\sO_{\nP(E_{k+1,L})}(1)|$ induces a morphism $\beta_k$ from $B^k(L):=\nP(E_{k+1,L})$ to $\nP H^0(C,L) = \nP^r$. If $\pi_k \colon B^k(L) \to C_{k+1}$ is the canonical projection, then $\beta_k(\pi_k^{-1}(\xi))$ is the $k$-plane in $\nP^r$ such that $\beta_k(\pi_k^{-1}(\xi)) \cap C = \xi$ for any $\xi \in C_{k+1}$. Thus $\beta_k(B^k(L))$ is the $k$-th secant variety $\Sigma_k=\Sigma_k(C, L)$ of $C \subseteq \nP^r$. We may think that $\beta_k \colon B^k(L) \to \Sigma_k$ is a resolution of singularities. Let $Z_{k-1}:=\beta_k^{-1}(\Sigma_{k-1})$. Then $Z_{k-1}= (k+1)H_k - \pi_k^* A_{k+1,L}$, where $\sO_{B^k(L)}(H_k) = \sO_{\nP(E_{k+1,L})}(1)$ and $A_{k+1,L}:=S_{k+1,L}(-2\delta_{k+1})$ (see \cite[Proposition 3.15]{ENP1}). Notice that $\omega_{B^k(L)}(Z_{k-1}) = \pi_k^* S_{k+1, \omega_C}$.
Ein--Niu--Park \cite{ENP1} established the following theorem (cf. Theorem \ref{thm:ENP}).

\begin{theorem}\label{thm:secant} Let $R:=R(\Sigma_k; \sO_{\Sigma_k}(1))$. We have the following:
\begin{enumerate}
\item[$(a)$] $R^i \beta_{k,*}\sO_{B^k(L)}(-Z_{k-1})=\begin{cases} \sI_{\Sigma_{k-1}|\Sigma_k} & \text{for $i=0$} \\ 0 & \text{for $i \geq 1$}. \end{cases}$
\item[$(b)$] $R^i \beta_{k,*}\omega_{B^k(L)}(Z_{k-1}) = \begin{cases} \omega_{\Sigma_k} & \text{for $i=0$} \\ \omega_{\Sigma_{k-1}} & \text{for $i=1$} \\ 0 & \text{for $i \geq 2$}. \end{cases}$
\item[$(c)$] $\Sigma_k$ has normal Du Bois singularities.
\item[$(d)$] $\Sigma_k \subseteq \nP^r$ is arithmetically Cohen--Macaulay and $H^{2k+1}(\Sigma_k, \sO_{\Sigma_k}(\ell))=0$ for $\ell>0$. In particular, $\reg(R)=2k+2$ unless $g=0$, in which case $\reg(R)=k+1$.
\item[$(e)$] If $\deg L \geq 2g+2k+1+p$, then $\sO_{\Sigma_k}(1)$ satisfies $N_{k+2,p}$-property.
\end{enumerate}
\end{theorem}

To prove the theorem, we proceed by the induction on $k$. As Ullery did in \cite{Ullery}, we apply the formal function theorem to prove the assertions $(a), (b), (c)$, which are local questions. Here we  focus on the normality of $\Sigma_k$. We remark that the idea of proving the normality of a variety using the formal function theorem was already presented in Ein's old paper \cite{Ein}. Note that $\Sigma_k$ is normal if and only if the map 
$\sO_{\nP^r} \to \beta_{k, *} \sO_{B^k(L)}$ is surjective at any point $x \in \Sigma_m \setminus \Sigma_{m-1}$ for $0 \leq m \leq k-1$. Notice that the unique $(m+1)$-secant plane containing $x \in \Sigma_m$ determines $\xi \in C_{m+1}$ uniquely. By the formal function theorem, we need to confirm that the map 
$$
\varprojlim \sO_{\nP^r}/\mathfrak{m}_x^{\ell} \longrightarrow \varprojlim H^0(\sO_{B^k(L)}/\sI_{F_x}^{\ell})
$$ 
is surjective, where $\mathfrak{m}_x$ is the ideal sheaf of $x$ in $\nP^r$ and $\sI_{F_x}$ is the ideal sheaf of $F_x:=\beta_k^{-1}(x)= C_{k-m}$ in $B^k(L)$. To this end, consider the commutative diagram
$$
\xymatrixcolsep{0.2in}
\xymatrix{
0\ar[r]&	\mathfrak{m}^{\ell}/\mathfrak{m}_x^{\ell+1} \ar[d]^-{\alpha_{\ell}} \ar[r] & \sO_{\nP^r}/\mathfrak{m}_x^{\ell+1} \ar[r]\ar[d]  & \sO_{\nP^r}/\mathfrak{m}_x^{\ell}\ar[d] \ar[r] & 0\\
0\ar[r]&	H^0(\sI_{F_x}^{\ell}/\sI_{F_x}^{\ell+1}) \ar[r]& H^0(\sO_{B^k(L)}/(\sI_{F_x}^{\ell+1})\ar[r]  & H^0(\sO_{B^k(L)}/\sI_{F_x}^{\ell})\ar[r] &\cdots
}
$$
and proceeding by induction on $\ell$, we reduce the problem to check that $\alpha_{\ell}$ is surjective for all $\ell \geq 1$. Note that $\mathfrak{m}_x^{\ell}/\mathfrak{m}_x^{\ell+1} = S^\ell T_x^\vee \nP^r$ and $\sI_{F_x}^{\ell}/\sI_{F_x}^{\ell+1} = S^{\ell} N_{F_x|B^k(L)}^{\vee}$
As $\alpha_1$ is an isomorphism (\cite[Proposition 3.13]{ENP1}) and $\alpha_{\ell} = \theta_{\ell} \circ S^{\ell} \alpha_1$, where 
$$
\theta_{\ell} \colon S^{\ell} H^0(N_{F_x|B^k(L)}^{\vee}) \longrightarrow H^0(S^{\ell} N_{F_x|B^k(L)}^{\vee}),
$$ 
it is sufficient to show that $\theta_{\ell}$  is surjective for $\ell \geq 1$. Based on Bertram's work \cite{Bertram}, it was shown in \cite[Proposition 3.13]{ENP1} that 
$$
N_{F_x|B^k(L)}^{\vee} = \sO_{F_x}^{\oplus 2m+1} \oplus E_{k-m, L(-2\xi)}.
$$
Then the surjectivity of $\theta_{\ell}$ easily follows from the projective normality of lower secant varieties $\Sigma_{k-m-1}(C, L(-2\xi)) \subseteq \nP H^0(C, L(-2\xi))$ (notice that $\deg L(-2\xi) \geq 2g + 2(k-m-1)+1$). Therefore, $\Sigma_k$ has normal singularities. The other assertions in $(a), (b), (c)$ can be similarly proven except two cases. First, after establishing all the other assertions, we can prove $\beta_{k,*}\omega_{B^k(L)}(Z_{k-1}) = \omega_{\Sigma_k}$ in $(b)$ using \cite[Theorem 1.1]{KSS} (see \cite[Proposition 3.6]{CKP}). Second, we apply a theorem of Koll\'{a}r--Kov\'{a}cs \cite[Corollary 6.28]{Kollar} and $(a)$ to show that $\Sigma_k$ has Du Bois singularities as claimed in $(c)$. For this purpose, we need to check that $Z_{k-1}$ has Du Bois singularities. In fact, one can show that $(B^k(L), Z_{k-1})$ is a log canonical pair (\cite[Proposition 3.15]{ENP1}) based on Bertram's blow-up construction \cite{Bertram} (see also \cite[Subsection 3.3]{ENP1}). This implies that $Z_{k-1}$ has Du Bois singularities by \cite[Corollary 6.32]{Kollar}.

Now, we prove the remaining assertions $(d),(e)$ assuming $(a),(b),(c)$. By Danila's theorem (see \cite{ENP2}), $H^0(\nP^r, \sO_{\nP^r}(\ell)) = H^0(\Sigma_k, \sO_{\Sigma_k}(\ell))$ for $1 \leq \ell \leq k+1$. Then one can easily see that
$$
H^i(\Sigma_k, \wedge^j M_{\sO_{\Sigma_k}(1)} \otimes \sO_{\Sigma_k}(k+1))=0~~\text{ for $i \geq 1$, $j \geq 0$, and $i \geq j-p$}
$$
implies that $(i)$ $H^i(\Sigma_k, \sO_{\Sigma_k}(\ell))=0$ for $i>0, \ell>0$, $(ii)$ $\Sigma_k \subseteq \nP^r$ is projectively normal, and $(iii)$ $\sO_{\Sigma_k}(1)$ satisfies $N_{k+2,p}$-property when $\deg L \geq 2g+2k+1+p$. Considering the short exact sequence 
$$
0 \longrightarrow \wedge^j M_{\sO_{\Sigma_k}(1)} \otimes \sI_{\Sigma_{k-1}|\Sigma_k} \longrightarrow \wedge^j M_{\sO_{\Sigma_k}(1)} \longrightarrow \wedge^j M_{\sO_{\Sigma_{k-1}}(1)} \longrightarrow 0
$$
and using the induction hypothesis, we transfer the problem to 
$$
H^i(\Sigma_k, \wedge^j M_{\sO_{\Sigma_k}(1)} \otimes \sI_{\Sigma_k}(k+1))=0~~\text{ for $i \geq 1$, $j \geq 0$, and $i \geq j-p$}.
$$
By the Du Bois-type condition $(a)$, it is equivalent to 
$$
H^i(B^k(L), \wedge^j M_{H_k} \otimes \pi_k^* A_{k+1, L})=0~~\text{ for $i \geq 1$, $j \geq 0$, and $i \geq j-p$},
$$
which can be identified with
$$
H^i(C_{k+1}, \wedge^j M_{k+1, L} \otimes A_{k+1, L})=0~~\text{ for $i \geq 1$, $j \geq 0$, and $i \geq j-p$},
$$
since $\pi_{k,*}\wedge^j M_{H_k} = \wedge^j M_{k+1, L}$ and $R^i \pi_{k,*}\wedge^j M_{H_k} = 0$ for $i>0$ (\cite[Lemma 5.1]{ENP1}). Let $\Delta_{k+1}:=q_{k+1}^*\delta_{k+1}$ be the union of all pairwise diagonals in $C^{k+1}$. As $\sO_{C_{k+1}}(-\delta_{k+1})$ is a direct summand of $q_{k+1,*}\sO_{C^{k+1}}$ (\cite[Lemma 3.5]{ENP1}), the desired vanishing follows from
$$
H^i(C^{k+1}, \wedge^j q_{k+1}^* M_{k+1,L} \otimes L^{\boxtimes k+1}(-\Delta_{k+1}))=0~~\text{ for $i \geq 1$, $j \geq 0$, and $i \geq j-p$}.
$$
Using Rathmann's technique \cite{Rathmann1}, we can prove this vanishing (see \cite[Theorem 4.1]{ENP1}). It only remains to show that $\Sigma_k \subseteq \nP^r$ is arithmetically Cohen--Macaulay. It is enough to verify that
\begin{equation}\label{eq:ACM}
H^i(\Sigma_k, \sO_{\Sigma_k}(-\ell))=0~~\text{ for $1 \leq i \leq 2k$ and $\ell \geq 0$}.
\end{equation}
The proof of (\ref{eq:ACM}) in \cite{ENP1} is incomplete as was pointed out by Doyoung Choi (we would like to thank him). By \cite[Theorem 10.42]{Kollar}, it was reduced to the cases of $\ell = 0,1$, and these cases are verified in the proof of \cite[Theorem 5.8]{ENP1}. Here \cite[Theorem 10.42]{Kollar} states that if $X$ is a projective pure dimensional scheme with Du Bois singularities and $A$ is an ample line bundle on $X$, then $h^i(X, A^{-\ell}) = h^i(X, A^{-1})$ for $i< \dim X$ and $\ell \geq 1$. This statement  contains a crucial typo: $h^i(X, A^{-\ell}) = h^i(X, A^{-1})$ should be replaced by $h^i(X, A^{-\ell}) \geq h^i(X, A^{-1})$. Thus  (\ref{eq:ACM}) was only proven for $\ell=0,1$. However, using the fact that 
$$
H^i(C_{k+1}, S_{k+1, \omega_C})=\begin{cases} S^{k+1} H^0(C, \omega_C) & \text{for $i=0$} \\ S^k H^0(C, \omega_C) & \text{for $i=1$} \\ 0 & \text{for $i \geq 2$} \end{cases}
$$
and the assertion $(b)$ (the former fact is also used to prove $(b)$), one can establish (\ref{eq:ACM}) for $\ell=0$ and $\ell \gg 0$ by similar arguments in the proof of \cite[Theorem 5.8]{ENP2}. Then the corrected version of \cite[Theorem 10.42]{Kollar} implies that (\ref{eq:ACM}) holds for all $\ell \geq 0$.

As explained in the introduction, Theorem \ref{thm:secant} determines the first $p+1$ columns of the Betti table of $R(\Sigma_k; \sO_{\Sigma_k}(1))$ when $\deg L = 2g+2k+1+p$. Assuming $\deg L \gg 0$, Choe--Kwak--Park \cite{CKP} determined the remaining part of the Betti table using the \emph{gonality sequence} $(\gamma^0(C)+0, \gamma^1(C)+1, \gamma^2(C)+2, \ldots)$ of $C$, where
$$
\gamma^q(C):=\min\{d -q \mid \text{$C$ carries a linear series $g_d^q$} \}.
$$
Note that $\gamma^1(C)+1=\gon(C)$.
 Let $e:=\codim \Sigma_k=r-2k-1$. Recall Theorem \ref{thm:CKP}: for each $k+1 \leq q \leq 2k+2$, if $e-g+1 \leq p \leq e$, then 
$$
K_{p,q}(\Sigma_k; \sO_{\Sigma_k}(1)) \neq 0~~\Longleftrightarrow~~e - g+1 \leq p \leq e - \gamma^{2k+2-q}(C).
$$
We give a sketch of the proof of this theorem. Using the Du Bois-type condition $(a)$ and proceeding by induction on $q-k-1$, we reduce the vanishing part to
$$
H^{q^*+1}(C_{k+2}, \wedge^{p^*+q^*} M_{E_{k+2, L}} \otimes S_{k+2, \omega_C})=0,
$$
where $p^*:=e-p$ and $q^*:=2k+2-q$ are introduced by duality theorem (\cite[Theorem 2.c.1]{Green1}).
By \cite[Lemma 3.3]{CKP}, it is equivalent to
$$
H^{q^*+1}(C_{p^*+q^*} \times C_{k+2}, (N_{p^*+q^*,L} \boxtimes S_{k+2, \omega_C})(-D_{p^*+q^*, k+2})) = 0.
$$
Let $\pr_1 \colon C_{p^*+q^*} \times C_{k+2} \to C_{p^*+q^*}$ be the projection map. By the Leray spectral sequence for $\pr_1$, it is enough to check that
$$
H^i(C_{p^*+q^*}, R^{q^*+1-i} \pr_{1,*} (N_{p^*+q^*,L} \boxtimes S_{k+2, \omega_C})(-D_{p^*+q^*, k+2})) =0~~
\text{ for $0 \leq i \leq q^*+1$}.
$$
When $i>0$, the cohomology vanishing follows from Fujita--Serre vanishing since $N_{p^*+q^*, L}$ is sufficiently positive. When $i=0$, we consider the fiber of $R^{q^*+1} \pr_{1,*} (N_{p^*+q^*,L} \boxtimes S_{k+2, \omega_C})(-D_{p^*+q^*,k+2})$ over $\xi \in C_{p^*+q^*}$, which is 
$$
H^{q^*+1}(C_{k+2}, S_{k+2, \omega_C(-\xi)}) = S^{k+1-q^*} H^0(C, \omega_C(-\xi)) \otimes \wedge^{q^*+1} H^1(C, \omega_C(-\xi)).
$$
Now, the gonality sequence condition $p \geq e-\gamma^{2k+2-q}(C)+1$ can be restated as $\gamma^{q^*}(C) \geq p^*+1$, which gives $h^1(C, \omega_C(-\xi)) \leq q^*$ in view of the positivity of $\omega_C$. Thus we obtain
$$
R^{q^*+1} \pr_{1,*} (N_{p^*+q^*,L} \boxtimes S_{k+2, \omega_C})(-D_{p^*+q^*,k+2}) = 0,
$$
which implies the desired cohomology vanishing for $i=0$. 

For the nonvanishing part, it suffices to show that the map
$$
\begin{array}{l}
H^{q-1}(\wedge^{p+q-1}M_{\sO_{\Sigma_k}(1)}\otimes \sI_{\Sigma_{k-1}|\Sigma_k}(1)) 
\longrightarrow H^q(\wedge^{p+q-1}M_{\sO_{\Sigma_{k+1}}(1)}\otimes \sI_{\Sigma_k|\Sigma_{k+1}}(1))
\end{array}
$$
is nonzero. We can reduce to showing that the map
$$
\begin{array}{l}
R^{q^*+1} \pr_{1,*} (N_{p^*+q^*,L} \boxtimes S_{k+2, \omega_C})(-D_{p^*+q^*,k+2}) \\
 \longrightarrow R^{q^*} \pr_{1,*} (N_{p^*+q^*,L} \boxtimes S_{k+1, \omega_C})(-D_{p^*+q^*,k+1}) \otimes H^1(C, \omega_C)
 \end{array}
$$
is nonzero. By the gonality sequence condition $p \leq e-\gamma^{2k+2-q}(C)$, which can be rephrased as $\gamma^{q^*}(C) \leq p^*$, we can find an effective divisor $\xi \in C_{p^*+q^*}$ with $h^0(C, \omega_C(-\xi)) = g-p^*$ and $h^1(C, \omega_C(-\xi)) = q^*+1$. The above map over $\xi \in C_{p^* + q^*}$ is $\id_{S^{k+1-q^*} H^0(C, \omega_C(-\xi))} \otimes \delta$, where $\delta$ is the Koszul-like map
$$
\wedge^{q^*+1} H^1(C, \omega_C(-\xi)) \longrightarrow \wedge^{q^*} H^1(C, \omega_C(-\xi)) \otimes H^1(C, \omega_C).
$$
Since $\delta$ is clearly nonzero, it follows that the above map is nonzero.

In view of the effective gonality theorem (Theorem \ref{thm:effgon}), it is a natural problem to establish an effective lower bound for $\deg L$ such that the conclusion of Theorem \ref{thm:CKP} holds. Some partial answers and related results are discussed in \cite[Section 5]{CKP}.

\section{Syzygies of tangent developable surfaces and generic Green's conjecture}\label{sec:tangent}
For a moment, we assume that $\operatorname{char}(\kk) = 0$. Let $C$ be a smooth projective curve of genus $g \geq 3$. If $\Cliff(C) \geq 1$ (i.e., $C$ is nonhyperelliptic), then $\omega_C$ is very ample and M. Noether's theorem says that $C \subseteq \nP H^0(C, \omega_C)=\nP^{g-1}$ is projectively normal. Petri's theorem states that if $\Cliff(C) \geq 2$, then $I_{C|\nP^{g-1}}$  is generated by quadrics. In the 1980s, Green \cite[Conjecture 5.1]{Green1} formulated a very famous conjecture predicting that $\omega_C$ satisfies $N_{(\Cliff(C)-1)}$-property. By the duality theorem (\cite[Theorem 2.c.1]{Green1}), this completely determines the shape of the Betti table of $R(C; \omega_C)$. This conjecture is still widely open. However, in the 2000s, Voisin \cite{Voisin1,  Voisin2} resolved Green's conjecture for general curves. For this purpose, it suffices to exhibit one curve of each genus for which the assertion holds. If $S \subseteq \nP^g$ is a K3 surface of degree $2g-2$, then its general hyperplane section is a canonical curve $C \subseteq \nP^{g-1}$ and $K_{p,2}(S; \sO_S(1)) = K_{p,2}(C; \omega_C)$. Voisin used the geometry of the Hilbert schemes of K3 surfaces for generic Green's conjecture. Recently, Kemeny \cite{Kemeny} gave a simpler proof of Voisin's theorem for even genus case and a streamlined version of her arguments for odd genus case.  

Prior to Voisin's work,  in 1980s, O'Grady and Buchweitz--Schreyer independently observed that one could use the tangent developable surface $T \subseteq \nP^g$ of a rational normal curve of degree $g$ to solve Green's conjecture for general curves of genus $g$. Note that $T \subseteq \nP^g$ is arithmetically Cohen--Macaulay and it is a degeneration of a K3 surface. A general hyperplane section of $T$ in $\nP^g$ is a canonically embedded $g$-cuspidal rational curve $\overline{C} \subseteq \nP^{g-1}$ of degree $2g-2$, which is a degeneration of a general canonical curve $C \subseteq \nP^{g-1}$ with $\Cliff(C) = \lfloor (g-1)/2 \rfloor$. By the upper semicontinuity of graded Betti numbers, $K_{p,2}(\overline{C}; \sO_{\overline{C}}(1)) = K_{p, 2}(T; \sO_T(1)) = 0$ implies $K_{p,2}(C; \omega_C)=0$. For generic Green's conjecture, it suffices to prove that
$$
K_{p,2}(T; \sO_T(1))=0~~\text{ for $0 \leq p \leq \lfloor (g-3)/2 \rfloor$}.
$$
This was finally established by Aprodu--Farkas--Papadima--Raicu--Weyman \cite{AFPRW} when $\operatorname{char}(\kk) = 0$ or $\geq (g+2)/2$ (see Theorem \ref{thm:AFPRW}). This important result gives not only an alternative proof of generic Green's conjecture but also an extension to positive characteristic. This circle of ideas is surveyed in \cite{EL4}.

We explain the proof of Theorem \ref{thm:AFPRW} given in \cite{AFPRW}. We have a short exact sequence
$$
\xymatrix{
0 \ar[r] & \sO_T \ar[r] & \nu_* \sO_{\widetilde{T}} \ar[r] & \omega_{\nP^1} \ar[r] & 0,
}
$$
where $\nu \colon \widetilde{T} = \nP^1 \times \nP^1 \to T$ is a resolution of singularities. It is easy to see that 
$$
K_{p,2}(T; \sO_T(1)) = \coker\big(K_{p,1}(T, \nu_* \sO_{\widetilde{T}}; \sO_T(1)) \xrightarrow{~\gamma~} K_{p,1}(\nP^1, \omega_{\nP^1}; \sO_{\nP^1}(g)) \big).
$$
Let $U:=H^0(\nP^1, \sO_{\nP^1}(1))$, $V:=D^{p+2} U$, $W:=D^{2p+2} U$, and $q:=g-p-3$. Here $D^m U$ is the $m$-th divided power of $U$ (see e.g., \cite[Section 3]{AFPRW}, \cite[Section 2]{Park2}). The authors of \cite{AFPRW} devoted considerable effort to show that $\gamma$ arises as the composition
\begin{equation}\label{eq:gammaintro}
S^{q} V \otimes W \xrightarrow{\id_{S^{q}V} \otimes \Delta}  S^{q} V \otimes \wedge^2 V \xrightarrow{~\delta~}  \ker(S^{q+1}V \otimes V \xrightarrow{~\delta~} S^{q+2}V),
\end{equation}
where $\Delta$ is the co-Wahl map and $\delta$ is the Koszul differential. To achieve this, they established an explicit characteristic-free Hermite reciprocity for $\mathfrak{sl}_2$-representations, and they carried out complicated algebraic computations. Now, $K_{p,2}(T, \sO_T(1))$ is the degree $q$ piece of the \emph{Koszul module} $W_q(V, W)$, which is defined to be the homology of a complex
$$
S^q V \otimes W \xrightarrow{~\gamma~} S^{q+1}V \otimes V \xrightarrow{~\delta~} S^{q+2} V.
$$
Then it is enough to prove that 
\begin{equation}\label{eq:vanKosmod}
W_q(V, W) = 0~~\text{ for $q \geq p$}.
\end{equation}
The vanishing result (\ref{eq:vanKosmod}) was established in the previous work of the authors of \cite{AFPRW} in characteristic zero and extended to positive characteristics in \cite{AFPRW}.

Park \cite{Park2} gave a geometric approach to proving Theorem \ref{thm:AFPRW} utilizing the first secant variety $\Sigma:=\Sigma_1(\nP^1, L)$ of a rational normal curve $C$ of degree $g$, where $L:=\sO_{\nP^1}(g)$. We remark that Theorem \ref{thm:secant} holds in positive characteristics when $C=\nP^1$. Note that $Q:=\sigma(D_{1,1})$ is a smooth conic in $\nP^2$, where $\sigma:=\sigma_{1,1} \colon \nP^1 \times \nP^1 \to \nP^2$ is the universal family over $(\nP^1)_2=\nP^2$. Let $\widetilde{T}:=\nP(E|_Q)$, where $E:=E_{2, L}$. Then $\widetilde{T}$ is linearly equivalent to $K_{B} + Z$, where $B:=B^1(L)$ and $Z:=Z_0=\nP^1 \times \nP^1$, and $\widetilde{T}=\nP^1 \times \nP^1$ when $\operatorname{char}(\kk) = 0$ or $\geq (g+2)/2$. Then $\beta_1(\widetilde{T})=T$ is the tangent developable surface of $v_g(\nP^1) \subseteq \nP H^0(\nP^1, L)=\nP^g$, and $\beta_1|_{\widetilde{T}}$ coincides with the map $\nu$. One can easily check that $T$ is a Weil divisor on $\Sigma$ linearly equivalent to $-K_{\Sigma}$. Let $H:=H_1$ be the tautological divisor on $B=\nP (E)$. We can realize (\ref{eq:gammaintro}) as maps induced in cohomology of vector bundles on $B$:
$$
H^1(\widetilde{T}, \wedge^{p+2} M_H|_{\widetilde{T}}) \xrightarrow{~\alpha~} H^2(B, \wedge^{p+2} M_H \otimes \omega_B(Z)) \xrightarrow{~\delta~} H^2(Z, \omega_{\nP^1} \boxtimes \wedge^{p+2} M_L \otimes \omega_{\nP^1}).
$$
It is easy to check that $\alpha = \id_{S^qV} \otimes \Delta$, where $\Delta$ is the dual of the restriction map $H^0(\nP^2, \sO_{\nP^2}(p+1)) \rightarrow H^0(Q, \sO_Q(p+1))$. Put $M_E:=M_{2,L}$. The map $\delta$ is naturally factored as\\[-20pt]

\begin{scriptsize}
$$
\xymatrixrowsep{0.1in}
\xymatrixcolsep{0.3in}
\xymatrix{
H^2(\wedge^{p+2} M_H \otimes \omega_B(Z)) \ar@{=}[d] \ar[r]^-{\id_{S^q V} \otimes \iota} & H^2(\sigma^* \wedge^{p+2} M_E \otimes (\omega_{\nP^1} \boxtimes \omega_{\nP^1})) \ar@{=}[d] \ar[r]^-{m \otimes \id_V} & H^2(\omega_{\nP^1} \boxtimes \wedge^{p+2} M_L \otimes \omega_{\nP^1})\ar@{=}[d] \\
S^q V \otimes \wedge^2 V & S^q V \otimes V \otimes V & S^{q+1} V \otimes V,
}
$$
\end{scriptsize}

\noindent where $\iota$ is the canonical injection and $m$ is the multiplication map identified with
$$
H^1(\nP^{p+2} \times \nP^1, \sO_{\nP^{p+2}}(q) \boxtimes \omega_{\nP^1}(-p-2)) \longrightarrow H^1(\nP^{p+2} \times \nP^1, \sO_{\nP^{p+2}}(q+1) \boxtimes \omega_{\nP^1}).
$$
This gives a quick verification of the descriptions of the maps in (\ref{eq:gammaintro}). Next, regarding $V=H^0(\nP^{p+2}, \sO_{\nP^{p+2}}(1))$, we can give a direct proof of (\ref{eq:vanKosmod}) as follows.
We put $M_V:=M_{\sO_{\nP^{p+2}}(1)}$. Then (\ref{eq:gammaintro}) can be identified with
$$
W \otimes H^0(\nP^{p+2}, \sO_{\nP^{p+2}}(q)) \longrightarrow \wedge^2 V  \otimes H^0(\nP^{p+2}, \sO_{\nP^{p+2}}(q)) \longrightarrow H^0(\nP^{p+2}, M_V(q+1)).
$$
We can show that $W \otimes \sO_{\nP^{p+2}} \to M_V(1)$ is surjective and its kernel is $(p+1)$-regular. This implies (\ref{eq:vanKosmod}). For this last part, we use the characteristic assumption $\operatorname{char}(\kk) = 0$ or $\geq (g+2)/2$.
It is worth noting that the characteristic assumption cannot be improved. Indeed, if $2 \leq \operatorname{char}(\kk) \leq (g+1)/2$, then $K_{\lfloor (g-3)/2 \rfloor, 2}(T, \sO_T(1)) \neq 0$ (see \cite[Remark 5.17]{AFPRW}).

Concerning the characteristic of the base field, Raicu--Sam \cite[Theorem 1.5]{RS} obtained a sharp result: Green's conjecture holds for general curves of genus $g$ when $\operatorname{char}(\kk) = 0$ or $\operatorname{char}(\kk) \geq \lfloor(g-1)/2 \rfloor$. This confirms a conjecture of Eisenbud--Schreyer. 
It has long been known that generic Green's conjecture would follow from the canonical ribbon conjecture. A \emph{canonical ribbon} is a certain double structure on a rational normal curve, and it can be realized as a hyperplane section of a \emph{K3 carpet} $X=X(a,b) \subseteq \nP^{a+b+1}$ for integers $b \geq a \geq 1$, which is a unique double structure on a rational normal surface scroll $S(a,b)=\nP (\sO_{\nP^1}(a) \oplus \sO_{\nP^1}(b)) \subseteq \nP^{a+b+1}$ such that $\omega_X = \sO_X$ and $h^1(X, \sO_X)=0$. Then $X$ is a degeneration of a K3 surface of degree $2(a+b)$, and its hyperplane section (a canonical ribbon) is a degeneration of a general canonical curve $C$ of genus $a+b+1$ with $\Cliff(C) = a$. By extending algebraic arguments of \cite{AFPRW}, Raicu--Sam  \cite[Theorem 1.1]{RS} proved that $\sO_X(1)$ satisfies $N_{(a-1)}$-property when  $\operatorname{char}(\kk) = 0$ or $\operatorname{char}(\kk) \geq a$. For Eisenbud--Schreyer's conjecture, we only need to consider the case $a=\lfloor (g-1)/2 \rfloor$ and $b=\lfloor g/2 \rfloor$. In this case, we take a general double line $\ell \in |Q|=|\sO_{\nP^2}(2)|$. Then $\beta_1(\nP(E|_{\ell}))=X(\lfloor (g-1)/2 \rfloor, \lfloor g/2 \rfloor)  \subseteq \nP^g$ is a K3 carpet. Along the same way as the geometric proof of Theorem \ref{thm:AFPRW}, we can obtain a simple geometric proof of Raicu--Sam's result \cite[Theorem 1.1]{RS}. 

\begin{theorem}[{$\operatorname{char}(\kk) = 0$ \normalfont{or} $\geq \lfloor (g-1)/2 \rfloor$ \normalfont{and} $\operatorname{char}(\kk) \neq 2$}]
Consider a K3 carpet $X=X(\lfloor (g-1)/2 \rfloor, \lfloor g/2 \rfloor) \subseteq \nP^g$ with $g \geq 3$. Then $\sO_X(1)$ satisfies $N_{\lfloor (g-3)/2 \rfloor}$-property.
\end{theorem}

This gives a very quick proof of generic Green's conjecture and Eisenbud--Schreyer's conjecture. So far, five different proofs of generic Green's conjecture are known: $\circled{1}$ Voisin's first proof using K3 surfaces \cite{Voisin1, Voisin2}, $\circled{2}$ Kemeny's simplification of Voisin's proof \cite{Kemeny} (Rathmann \cite{Rathmann2} provided a different account of Kemeny’s proof for even genus case), $\circled{3}$ Aprodu--Farkas--Papadima--Raicu--Weyman's algebraic proof using tangent developable surfaces of rational normal curves \cite{AFPRW}, $\circled{4}$ Raicu--Sam's algebraic proof using K3 carpets \cite{RS}, $\circled{5}$ Park's geometric simplification of the main results of \cite{AFPRW} and \cite{RS}. We remark that $\circled{1},\circled{2}$ only work in characteristic zero while $\circled{3}, \circled{4}, \circled{5}$ also work in positive characteristics. All proofs use K3 or K3-like surfaces. It would be exceedingly interesting to find a direct proof through symmetric products of curves. This may lead to the resolution of Green's conjecture for all curves.

It is quite natural to study syzygies of tangent developable surfaces of arbitrary smooth projective curves of genus $g \geq 1$. As a first step, Park \cite[Theorem 1.3]{Park2} showed the arithmetic normality of tangent developable surfaces.

\begin{theorem}[{$\operatorname{char}(\kk)=0$}]\label{thm:arithnormtansurf}
Let $C$ be a smooth projective curve of genus $g \geq 1$, and $L$ be a line bundle on $C$ with $\deg L \geq 4g+3$. Consider the tangent developable surface $T$ of $C$ embedded in $\nP H^0(C, L)=\nP^r$. Then $T \subseteq \nP^r$ is arithmetically normal but not arithmetically Cohen--Macaulay, and $H^i(T, \sO_T(m))=0$ for $i>0, m>0$ but $H^1(T, \sO_T) \neq 0, H^2(T, \sO_T) \neq 0$.
\end{theorem}

To prove the theorem, regarding $T$ as a Weil divisor on the first secant variety of $C$ in $\nP^r$, we transfer the problem to cohomology vanishing of vector bundles on the second symmetric (or ordinary) product of $C$. The most difficult part is to check the $2$-normality of $T \subseteq \nP^r$, which is turned out to be equivalent to the cohomology vanishing $H^1(C \times C, (L \boxtimes L)(-3D_{1,1})) = 0$. This was previously established by Bertram--Ein--Lazarsfeld \cite{BEL} when $\deg L \geq 4g+3$. The degree condition in Theorem \ref{thm:arithnormtansurf} is optimal. If $\deg L = 4g+2$, then $T \subseteq \nP^r$ is arithmetically normal if and only if $C$ is neither elliptic nor hyperelliptic (see \cite[Remark 4.4]{Park2}). Some relevant conjectures and examples on syzygies of tangent developable surfaces are presented in \cite[Section 4]{Park2}. For example, it was conjectured that if $g=1$, then $K_{p,1}(T; \sO_T(1))=0$ for $p \geq \lfloor (\deg L - 3)/2\rfloor$ and $K_{p,2}(T; \sO_T(1))=0$ for $0 \leq p \leq \lfloor (\deg L - 7)/3 \rfloor$. In this case, $T$ is a degeneration of an abelian surface, so this conjecture may have some applications to syzygies of general abelian surfaces.

\section{Asymptotic syzygies of algebraic varieties}\label{sec:asymptotic}
In this section, we assume that $\operatorname{char}(\kk) \geq 0$. Let $X$ be an $n$-dimensional smooth projective variety, $B$ be a line bundle on $X$, and $L_d:=\sO_X(dA+P)$ for an integer $d \gg 0$, where $A$ is an ample divisor and $P$ is an arbitrary divisor on $X$. Put $R_d:=R(X, B; L_d)$ and $r_d:=h^0(X, L_d)-1$. Assume that $d \gg 0$. Then $r_d = \Theta(d^n)$, and $H^i(X, B \otimes L^m)=0$ for $i>0$ and $m>0$. Thus $\reg (R_d) \leq n+1$, so $K_{p,q}(X, B; L)=0$ for $q \geq n+2$. 
We want to determine vanishing and nonvanishing of $K_{p,q}(X, B; L_d)$ for $0 \leq q \leq n+1$. By (\ref{eq:q=0,n+1}) in the introduction, we only need to focus on the cases of $1 \leq q \leq n$. Motivated by the nonvanishing results of \cite{EGHP, OP}, Ein--Lazarsfeld \cite[Theorem 4.1]{EL2} established the following:

\begin{theorem}[{Asymptotic Nonvanishing Theorem}]\label{thm:ant}
For each $1 \leq q \leq n$, there are constants $C_1, C_2 > 0$ such that if $d$ is sufficiently large, then
$$
K_{p,q}(X, B; L_d) \neq 0~~\text{ for $C_1 d^{q-1} \leq p \leq r_d - C_2 d^{n-1}$}.
$$
\end{theorem}

The influential paper \cite{EL2} opens the door to research on the asymptotic behaviors of $K_{p,q}(X, B; L_d)$ for $d$ increasing. It is natural to ask whether $K_{p,q}(X, B; L_d)=0$ for the values of $p$ outside the range in Theorem \ref{thm:ant}. Based on Raicu's result \cite[Appendix]{Raicu}, Park \cite[Theorem 1.1]{Park1} proved the following theorem, which was previously conjectured by Ein--Lazarsfeld \cite[Conjecture 7.1]{EL2}. 

\begin{theorem}[{Asymptotic Vanishing Theorem}]\label{thm:avt}
For each $1 \leq q \leq n$, there is a constant $C_3 > 0$ such that if $d$ is sufficiently large, then
$$
K_{p,q}(X, B; L_d) \neq 0~~\text{ for $p \leq C_3 d^{q-1}$}.
$$
\end{theorem}

Park \cite{Park3} subsequently realized that the asymptotic vanishing theorem leads to a quick proof of the more precise version of the asymptotic nonvanishing theorem, that is Theorem \ref{thm:asysyz}: for each $1 \leq q \leq n$, there are functions $c_q(d)$ and $c_q'(d)$ with
\begin{equation}\label{eq:c_qc_q'}
c_q(d) = \Theta(d^{q-1})~~\text{ and }~~c_q'(d) = \begin{cases} \Theta(d^{n-q}) & \text{if $H^{q-1}(X, B) =0$ or $q=1$} \\  q-1 & \text{if $H^{q-1}(X, B) \neq 0$ and $q \geq 2$} \end{cases}
\end{equation}
such that if $d$ is sufficiently large, then
$$
K_{p,q}(X, B; L_d) \neq 0~~\Longleftrightarrow~~c_q(d) \leq p \leq r_d-c_q'(d).
$$
It is natural to wonder whether $c_q(d)$ and $c_q'(d)$ are polynomials in $d$ when $d \gg 0$ (see \cite[Question 5.5]{Park3}). It seems very challenging even for $X=\nP^n$ (cf. \cite[Conjecture 2.3]{EL3}). However, the answer is known for $q=1$ by Yang \cite[Theorem 1]{DavidYang}. Furthermore, Ein--Lazarsfeld--Yang \cite[Theorem A]{ELY} proved that if $B$ is $p$-jet very ample, then $c_1(d) \geq p+1$, and Agostini \cite[Theorem A]{Agostini2} proved that if $B$ is not $p$-very ample, then $c_1(d) \leq p+1$. These results are already highly nontrivial and in fact a higher dimensional generalization of the gonality conjecture. We refer to \cite[Section 5]{Park3} for more discussion on related open problems.

Following \cite{Park1, Park3}, we give an outline of the proof of Theorem \ref{thm:asysyz}. First, we prove the asymptotic vanishing theorem for products of projective spaces (see \cite[Subsection 3.2]{Park1}). Let $Y:=\nP^{n_1} \times \cdots \times \nP^{n_{k-1}}$ and $X:=Y \times \nP^{n_k}$. Set 
$$
\begin{array}{l}
B_Y:=\sO_{\nP^{n_1}}(b_1) \boxtimes \cdots \boxtimes \sO_{\nP^{n_{k-1}}}(b_{k-1}),~ L_Y:=\sO_{\nP^{n_1}}(d_1) \boxtimes \cdots \boxtimes \sO_{\nP^{n_{k-1}}}(d_{k-1});\\
B:=B_Y \boxtimes \sO_{\nP^{n_k}}(b_k),~L:=L_Y \boxtimes \sO_{\nP^{n_k}}(d_k).
\end{array}
$$ 
For simplicity, we put $n:=n_k$. Fix $2 \leq q \leq n$. Assuming $d:=\min\{d_1, \ldots, d_k\} \gg 0$ and $p \leq (1/n_1! \cdots n_k!)(d^{q-1}+bd^{q-2}) = \Theta(d^{q-1})$, we show that
$$
K_{p,q}(X, B;L) = H^{q-1}\big(Y \times \nP^n, \wedge^{p+q-1}M_{L} \otimes ((L_Y \otimes B_Y) \boxtimes \sO_{\nP^n}(d_k+b_k))\big) = 0.
$$
Regarding $\nP^n=(\nP^1)_n$, let $\sigma:=\id_Y \times \sigma_{n-1,1}$, where $\sigma_{n-1,1} \colon \nP^{n-1} \times \nP^1 \to \nP^n$ is the universal family. It suffices to show that\\[-20pt]

\begin{small}
$$
H^{q-1}\big(Y \times \nP^{n-1} \times \nP^1, \wedge^{p+q-1}\sigma^* M_{L} \otimes ((L_Y\otimes B_Y) \boxtimes \sO_{\nP^{n-1}}(d_k+b_k) \boxtimes \sO_{\nP^1}(d_k+b_k+n-1))\big)=0.
$$
\end{small}

\noindent To this end, using the geometry of symmetric products of $\nP^1$, Park \cite[Subsection 3.1]{Park1} constructed a short exact sequence
$$
0 \longrightarrow  \bigoplus \sO_{Y \times \nP^{n-1}} \boxtimes \sO_{\nP^1}(-n) \longrightarrow  \sigma^* M_{L} \longrightarrow  M_{L_Y \boxtimes \sO_{\nP^{n-1}}(d_k)} \boxtimes \sO_{\nP^1}(d_k) \longrightarrow  0.
$$
Considering the filtration for $\wedge^{p+q-1}\sigma^* M_{L}$, we can reduce the problem to
$$
H^{q-1}\big(Y \times \nP^{n-1} \times \nP^1, (\wedge^i M_{L_Y \boxtimes \sO_{\nP^{n-1}}(d_k)} \otimes ((L_Y \otimes B_Y) \boxtimes \sO_{\nP^{n-1}}(d_k+b_k))) \boxtimes \sO_{\nP^1}(a_i) \big)= 0
$$
for $i \leq p+q-1$, where $a_i:= i(d_k+n) + d_k + b_k +2n-qn -1 -pn$.
By the K\"{u}nneth formula, it is equivalent to 
$$
\begin{array}{l}
(i)~H^{q-1-j}\big(Y \times \nP^{n-1}, \wedge^i M_{L_Y \boxtimes \sO_{\nP^{n-1}}(d_k)} \otimes ((L_Y \otimes B_Y) \boxtimes \sO_{\nP^{n-1}}(d_k+b_k))\big) = 0~\text{ or }\\
(ii)~H^j(\nP^1,  \sO_{\nP^1}(a_i))=0~~\text{ for $0 \leq i \leq p+q-1$ and $j=0,1$}.
\end{array}
$$
Assume $3 \leq q \leq n$. We proceed by induction on $n_1 + \cdots + n_k$. When $j=0$, $(i)$ holds for $i \leq p+q-1$ by induction. When $j=1$, $(i)$ holds for $i \leq (1/n_1! \cdots n_{k-1}!(n_k-1)!)(d^{q-2}+bd^{q-3})+q-2=\Theta(d^{q-2})$ by induction, and $(ii)$ holds for $i \geq (1/n_1! \cdots n_{k-1}!(n_k-1)!)(d^{q-2}+bd^{q-3})+q-1$ since $a_i \geq -1$. When $q=2$, the argument is slightly different but easier (see \cite[Subsection 3.2]{Park1}).

Next, we deduce Theorem \ref{thm:avt} for general case from the same assertion for products of three projective spaces proven in the previous paragraph. This reduction argument was done by Raicu  \cite[Appendix]{Raicu} (see also \cite[Subsection 3.3]{Park1}) shortly after Ein--Lazarsfeld's asymptotic vanishing conjecture \cite[Conjecture 7.1]{EL2} was raised.
We can choose suitable positive integers $a_1, a_2, a_3$ depending only on $X,A,P$ and positive integers $d_1, d_2, d_3$ grow linearly in $d$ in such a way that  if 
$$
\begin{array}{l}
A_1:=\sO_X(a_1A), ~A_2:=\sO_X(a_2A + P), ~A_3:=\sO_X(a_3A-P),\\
\nP^{n_1}:=\nP H^0(X, A_1),~\nP^{n_2}:=\nP H^0(X, A_2),~\nP^{n_3}:=\nP H^0(X, A_3),\\
Y:=\nP^{n_1} \times \nP^{n_2} \times \nP^{n_3}, ~L_Y:=\sO_{\nP^{n_1}}(d_1) \boxtimes \sO_{\nP^{n_2}}(d_2) \boxtimes \sO_{\nP^{n_3}}(d_3),
\end{array}
$$
then $L_d=A_1^{d_1} \otimes A_2^{d_2} \otimes A_3^{d_3}$ and there is a commutative diagram
$$
\xymatrixcolsep{0.9in}
\xymatrix{
X \ar@{^{(}->}[r] \ar@{^{(}->}[d] & Y \ar@{^{(}->}[d]\\
\nP^r:=\nP H^0(X, L_d) \ar@{^{(}->}[r]_-{\text{linear subspace}} & \nP H^0(Y, L_Y)=:\nP^N.
}
$$
We may regard $B$ as a coherent sheaf on $Y, \nP^r,\nP^N$. 
Since
$$
\begin{array}{l}
\text{(syzygies of $B$ on $\nP^N$)} = \text{(syzygies of $B$ on $\nP^r$)} \otimes \text{(Koszul complex of linear forms)},
\end{array}
$$
it follows that
$$
\min\{ p \mid K_{p,q}(X, B; L_d) \neq 0 \} = \min\{ p \mid K_{p,q}(Y, B; L_Y) \neq 0\}.
$$
By considering the minimal free resolution of $B$ on $Y$ with $\Z^3$-grading, we can deduce Theorem \ref{thm:avt} (the case of $q=1$ is trivial) from the result proven in the previous paragraph. Actually this argument shows that Theorem \ref{thm:avt} holds when $X$ is a possibly singular projective variety and $B$ is an arbitrary coherent sheaf on $X$.

Finally, we prove Theorem \ref{thm:asysyz} using the asymptotic vanishing theorem (Theorem \ref{thm:avt}). As in \cite{EL2}, the argument proceeds by induction on $n=\dim X$. Let $H$ be a suitably positive very ample line bundle on $X$, and choose a general member $\overline{X} \in |H|$. Put
$$
\overline{L}_d:=L_d|_{\overline{X}},~~\overline{B}:=B|_{\overline{X}},~~\overline{H}:=H|_{\overline{X}}, ~~V_d':=H^0(X, L_d \otimes H^{-1}).
$$
We have noncanonical splitting
$$
\begin{array}{l}
K_{p,q}(X, \overline{B}; L_d) = \bigoplus_{j=0}^p \wedge^j V_d' \otimes K_{p-j, q}(\overline{X}, \overline{B}; \overline{L}_d);\\
K_{p,q}(X, \overline{B} \otimes \overline{H}; L_d) =   \bigoplus_{j=0}^p \wedge^j V_d' \otimes K_{p-j, q}(\overline{X}, \overline{B} \otimes \overline{H}; \overline{L}_d).
\end{array}
$$
There are exact sequences\\[-22pt]

\begin{scriptsize}
$$
\begin{array}{l}
K_{p+1, q-1}(X, B \otimes H; L_d) \rightarrow K_{p+1, q-1}(X, \overline{B} \otimes \overline{H}; L_d) \xrightarrow{\theta_{p,q}} K_{p,q}(X, B; L_d) \rightarrow K_{p,q}(X, B \otimes H; L_d);\\[3pt]
K_{p,q}(X, B \otimes H^{-1}; L_d) \rightarrow K_{p,q}(X, B; L_d) \xrightarrow{\theta_{p,q}'} K_{p,q}(X, \overline{B}; L_d) \rightarrow K_{p-1, q+1}(X, B \otimes H^{-1}; L_d).
\end{array}
$$
\end{scriptsize}

\noindent When $q=1$ or $2$, the map $\theta_{p,q}$ should be slightly modified (see \cite[Section 3]{Park3}). By Theorem \ref{thm:avt} and induction on $n$, we can easily get the estimation of $c_q(d)$ and $c_q'(d)$ in (\ref{eq:c_qc_q'}). Furthermore, $\theta_{c_q(d), q}$ and $\theta_{r_d-c_q'(d), q}$ are nonzero maps. Thus there are\\[-20pt]

\begin{footnotesize}
$$
\begin{array}{l}
\text{$\alpha \in K_{c_q(d)+1-j_0, q-1}(\overline{X}, \overline{B} \otimes \overline{H}; \overline{L}_d)$ for some $0 \leq j_0 \leq c_q(d)+1$ lifted to $K_{c_q(d),q}(X, B; L_d)$;}\\
\text{ $\beta \in K_{r_d-c_q'(d)-j_0', q}(\overline{X}, \overline{B}; \overline{L}_d)$ for some $r_d' - c_q'(d) \leq j_0' \leq r_d'$ lifted to $K_{r_d-c_q'(d),q}(X, B; L_d)$.}
\end{array}
$$
\end{footnotesize}

\noindent Then we can lift $\alpha$ and $\beta$ to $K_{p,q}(X, B; L_d)$ for $c_q(d) \leq p \leq r_d - \Theta(d^{n-1})$ and $\Theta(d^{n-1}) \leq p \leq r_d-c_q'(d)$, respectively (see \cite[Theorem 3.1]{Park3} for more details on lifting syzygies from hypersurfaces). This gives the desired nonvanishing for $K_{p,q}(X, B; L_d)$ as stated in Theorem \ref{thm:asysyz}.

We close this paper by recalling open problems on the quantitative behavior of the nonzero graded Betti numbers $\kappa_{p,q}(X, B; L_d)$ as $d \to \infty$. Ein--Erman--Lazarsfeld \cite[Conjecture B]{EEL} (resp. Park \cite[Conjecture 5.7]{Park3}) conjectured that the Betti numbers $\kappa_{p,q}(X,B;L_d)$ are normally distributed (resp. form a unimodal sequence) for each $1 \leq q \leq n$. It is also expected that the Betti numbers $\kappa_{p,0}(X,B;L_d)$ are increasing and the Betti numbers $\kappa_{p,n+1}(X,B;L_d)$ are decreasing (see \cite[Remark 5.9]{Park3}). These have been verified when $n=1$ (see \cite[Proposition A]{EEL} and \cite[Proposition 5.8]{Park3}) but they are widely open when $n \geq 2$. It seems that these problems are already very challenging even for $X=\nP^n$. However, from the perspective of Boij--S\"{o}derberg theory \cite{BS, ES}, Ein--Erman--Lazarsfeld \cite[Corollary D]{EEL} proved that each row of the ``random'' Betti table is normally distributed.

\begin{acknowledgement}
The author would like to thank Daniele Agostini, Junho Choe, Lawrence Ein, Sijong Kwak, Wenbo Niu for the many interesting discussions we have had on syzygies of algebraic varieties. The author was partially supported by the National Research Foundation (NRF) funded by the Korea government (MSIT) (NRF-2021R1C1C1005479 and NRF-2022M3C1C8094326).\\[-15pt]
\end{acknowledgement}

\ethics{Competing Interests}{None}


\end{document}